\documentclass{mathincs}
\usepackage[T1]{fontenc}
\usepackage{floatflt}
\usepackage{algorithm,algpseudocode}
\usepackage{graphicx,caption}
\usepackage{amsfonts}
\usepackage[unicode, pdftex]{hyperref}
\usepackage{array}
\usepackage{longtable}
\usepackage{multirow}
\usepackage{vcell}
\usepackage[all]{xy}

 \newtheorem{thm}{Theorem}[section]
 \newtheorem{theorem}[thm]{Theorem}
 \newtheorem{example}[thm]{Example}

 \newtheorem{conjecture}[thm]{Conjecture}
 \theoremstyle{definition}
 \newtheorem{definition}[thm]{Definition}
 \theoremstyle{remark}

 \numberwithin{equation}{section}

\begin{document}
\title{An introduction to computational aspects of polynomial amoebas -- a survey}
%
%
\author{Vitaly A. Krasikov}

\address{%
Plekhanov Russian University of Economics \\
36 Stremyanny Lane\\
Moscow, Russia}
%
\email{krasikov.va@rea.ru,vitkras@inbox.ru}
%
\thanks{The research was supported by a grant from the Russian Science Foundation No.22-21-00556, \url{https://rscf.ru/project/22-21-00556/}}


\maketitle              

\begin{abstract}
This article is a survey on the topic of polynomial amoebas. We review results of papers written on the topic with an emphasis on its computational aspects. Polynomial amoebas have numerous applications in various domains of mathematics and physics. Computation of the amoeba for a~given polynomial and describing its properties is in general a problem of high complexity. We overview existing algorithms for computing and depicting amoebas and geometrical objects associated with them, such as contours and spines. We review the latest software packages for computing polynomial amoebas and compare their functionality and performance.

\end{abstract}

\section{Introduction}

The amoeba of a polynomial is the projection of its zero locus onto the space of absolute values in the logarithmic scale. Initially this notion was introduced by I.M. Gelfand, M.M. Kapranov, and A.V.~Zelevinsky in 1994~\cite{GKZ:1994}. Polynomial amoebas have numerous applications in topology~\cite{Guilloux-Marche:2021}, dynamical systems~\cite{Theobald:2002}, algebraic geometry~\cite{Goucha-Gouveia:2021,Lang:2019,Mikhalkin:2017}, complex analysis~\cite{Forsberg-Passare-Tsikh:2000}, mirror symmetry~\cite{Feng-He-Kennaway-Vafa:2008,Hicks:2020,Ruan:2000}, measure theory~\cite{Mikhalkin-Rullgard:2001,Passare-Rullgard:2000}, and statistical physics~\cite{Maslov:2015-1,Passare-Pochekutov-Tsikh:2013}.

There are several surveys published on the topic of polynomial amoebas (see~\cite{Mikhalkin:2004,deWolff:2017,Yger:2012}), but none of them overviewing the algorithms and the software for computing amoebas. The present paper is an attempt to make such an overview focused mainly on computational aspects of polynomial amoebas. The computation of amoebas is a~problem of great importance, and it also has a~high complexity for the following reasons. First, the amoeba of a~polynomial is an unbounded subset in~$\mathbb{R}^n$ (except for one-dimensional case) but at the same time in higher dimensions the amoeba for any given polynomial is a~very ``narrow'' subset of~$\mathbb{R}^n,$ so even locating the amoeba in the real space is a~non-trivial problem. The problem of deciding whether a~given point belongs to the amoeba (the membership problem) in general cannot be solved in polynomial time~\cite{Avendano-etAl:2018}. Since the membership problem is crucial for computing amoebas, there are numerous papers containing methods for solving it~\cite{Forsgard-Matusevich-Mehlhop-deWolff:2018,Purbhoo:2008,Theobald:2002,Theobald-deWolff:2015}). Another challenging problem is locating connected components of the amoeba complement, since their number can vary and they can be arbitrarily small. This is the reason why the computation of polynomial amoebas in general requires very high precision~\cite{Zhukov-Sadykov:2023}.

There is a~strong connection between the geometry of the amoeba of a~polynomial and its Newton polytope. Recall that the Newton polytope of a polynomial is the convex hull in~$\mathbb{R}^n$ of its support (that is, the set of its monomial exponent vectors). It is proved in~\cite{Forsberg-Passare-Tsikh:2000} that for any polynomial the number of connected components of its amoeba complement is at least equal to the number of vertices of its Newton polytope and at most equal to the total number of integer points contained in its Newton polytope. Amoebas with the minimal number of connected components in their complement are called~{\it solid} and amoebas with the maximal number of components are said to be {\it optimal.} 

There are several open problems related to polynomial amoebas, one of them considered in the present paper is the following conjecture by M. Passare~\cite{PassareConjectureFormulation:2009}: 

\begin{conjecture} (The Passare conjecture). \rm
\label{conj:Passare}
Let $p(z)$ be a maximally sparse polynomial (that is, the support of $p$ is equal to the set of vertices of its Newton polytope).  Then the amoeba of $p(z)$ is solid.
\end{conjecture}

The main purpose of the present survey is to overview and discuss the concept of a polynomial amoeba in the computational context, including historical details on the origin of polynomial amoebas, their modern applications, and computational aspects. It includes a review of algorithms and software packages for computing amoebas and testing these packages against different classes of polynomials. The main question we want to answer in the context of the Passare conjecture is whether the state of the art of the software for computing polynomial amoebas provides the means necessary to help in proving the conjecture or refuting it by finding a~counterexample. 

\section{Notation, basic properties, and examples}
\label{sec:basic_notation}

Let $p (z)$ be a (Laurent) polynomial in $n$ complex variables:
$$p(z_1,\ldots, z_n)=\sum\limits_{\alpha\in A} c_{\alpha} z^{\alpha}=\sum\limits_{\alpha\in A} c_{\alpha_1\ldots\alpha_n}z_1^{\alpha_1}\cdot\ldots\cdot z_n^{\alpha_n},$$
where $A \subset \mathbb{Z}^n$ is a finite set. Hereinafter the set~$A$ is called {\it the support} of $p(z).$ 

\begin{definition} \rm
{\it The amoeba} $\mathcal{A}_{p}$ of a polynomial $p(z)$ is the image of its zero locus under the map $\mathrm{Log}: \left(\mathbb{C^*}\right)^n \rightarrow \mathbb{R}^n,$ where $\mathbb{C^*}=\mathbb{C}\backslash \{0\}:$ $$\mathrm{Log}: (z_1, \ldots, z_n) \longmapsto (\mathrm{ln}|z_1|, \ldots, \mathrm{ln}|z_n|).$$
\end{definition}

Polynomial amoebas were introduced by I.M. Gelfand, M.M. Kapranov, and A.V. Zelevinsky in~\cite{GKZ:1994}. The reason for using the term ``amoeba'' is the resemblance between the shapes of the image $\mathrm{Log}({z|p(z) = 0})$ for $n = 2$ and the unicellular organism of the same name, in particular, a polynomial amoeba also has straight narrowing ``tentacles'' reaching to infinity. By using the $\mathrm{Log}$ map the notion of amoeba can be applied to algebraic or even transcendental hypersurfaces~\cite{Passare-Pochekutov-Tsikh:2013}, but in what follows we focus on the polynomial amoebas. 

The structure and the properties of the amoeba of a polynomial are closely connected to the structure of its Newton polytope.

\begin{definition} \rm
The convex hull in~$\mathbb{R}^n$ of the set~$A$ is called {\it the Newton polytope of} $p(z)$ (or {\it the Newton polygon of} $p(z)$ for $n=2$). Hereinafter we denote it by~$\mathcal{N}_p.$ 
\end{definition}

The connected components of the amoeba complement $\textrm{}^c\mathcal{A}_p = \mathbb{R}^n\backslash \mathcal{A}_p$ are convex subsets in~$\mathbb{R}^n.$ They are in bijective correspondence with the different Laurent expansions (centered at the origin) of the rational function~$1/p(z)$~\cite{GKZ:1994}. In the general case, the following statement holds:

\begin{theorem} \rm~\cite{Forsberg-Passare-Tsikh:2000}
The number of connected components of the amoeba complement $\textrm{}^c\mathcal{A}_p$ is at least equal to the number of vertices of the Newton polytope $\mathcal{N}_p$ and at most equal to the total number of integer points in $\mathcal{N}_p\cap \mathbb{Z}^n.$
\label{thm:connectedComponentsNumber}
\end{theorem}

The part of Theorem~\ref{thm:connectedComponentsNumber} concerning the lower bound for the number of connected components was proved in~\cite{GKZ:1994} and~\cite{Mkrtchian-Yuzhakov:1982}.

\begin{definition} \rm
If the number of connected components of~$\textrm{}^c\mathcal{A}_p$ is equal to the number of vertices of the Newton polytope~$\mathcal{N}_p,$ the amoeba~$\mathcal{A}_p$ is called {\it solid.} If the number of connected components of~$\textrm{}^c\mathcal{A}_p$ is equal to the number of integer points in~$\mathcal{N}_p\cap \mathbb{Z}^n,$ the amoeba~$\mathcal{A}_p$ is called {\it optimal.}
\end{definition}

\begin{example} \rm
Consider the family of polynomials $V_\alpha$ defined through the generating function (see~\cite[Chapter 2.3]{Xu-Dunkl:2014}) $$(1-2\langle b,x\rangle+||b||^2)^\frac{n-1}{2}=\sum\limits_{\alpha\in\mathbb{N}^n_0}b^\alpha V_\alpha(x),$$ where $b\in \mathbb{R}^n, ||b||<1, \langle\cdot,\cdot\rangle$ is the standard scalar product in~$\mathbb{R}^n.$

By neglecting a monomial factor of $V_\alpha$ and replacing variables $x_j, j=1,\ldots,n$ by $\xi_j=x_j^2,$ one obtains the polynomial~$\tilde V_\alpha(\xi)$ with the same number of components in its amoeba complement as $V(x)$ (see Lemma 2.6 and Example 2.10 in~\cite{Bogdanov-Sadykov:2020}).

The Newton polytope and the amoeba of the polynomial $\tilde V_{(5,2,3)}(\xi)$ are shown in Figure~\ref{fig:amoeba_3d}. Components of the amoeba complement are depicted there as colored convex shapes, the amoeba itself is a white space in between. The Newton polytope for~$\tilde V_{(5,2,3)}(\xi)$ contains~$12$ integer points and this number coincides with the number of connected components in the amoeba complement, so the amoeba of the polynomial $\tilde V_{(5,2,3)}(\xi)$ is optimal.

\begin{figure}[h!]
\begin{minipage}{6.5cm}
\includegraphics[width=6.5cm]{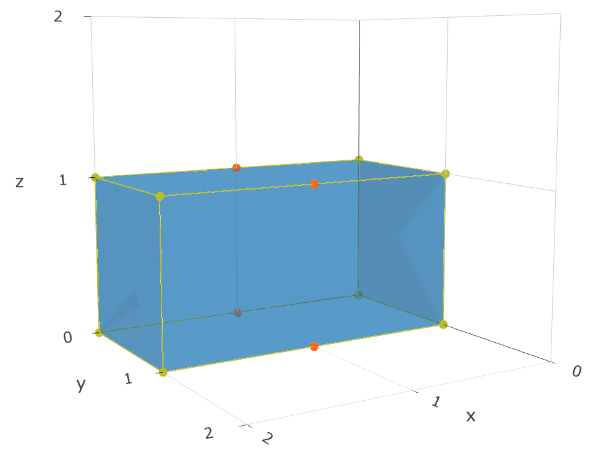}
\centering $(a)$
\end{minipage}
\begin{minipage}{7cm}
\includegraphics[width=7cm]{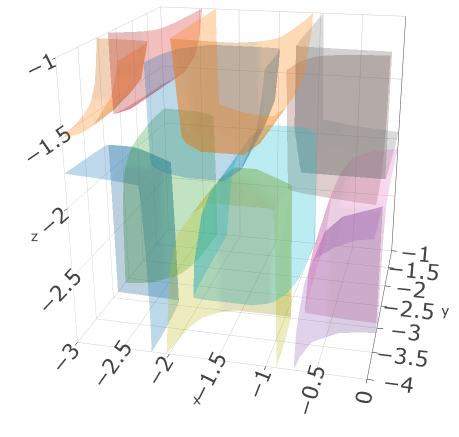}
\centering $(b)$
\end{minipage}

\caption{The Newton polytope (a) and connected components of the complement to the amoeba (b) of the polynomial~$\tilde V_{(5,2,3)}$}
\label{fig:amoeba_3d}
\end{figure}

\end{example}

It follows from Theorem~\ref{thm:connectedComponentsNumber} that there exists an injective function from the set of connected components of $\textrm{}^c\mathcal{A}_p$ into $\mathcal{N}_p\cap \mathbb{Z}^n.$ One may construct such a function by using the Ronkin function $$N_p(x)=\left(\frac{1}{2\pi i}\right)^n \int_{\mathrm{Log}^{-1}(x)}\mathrm{ln}|p(z)|\frac{dz}{z},\quad x\in\mathbb{R}^n.$$

The function $N_p$ is convex on $\mathbb{R}^n$ and it is affine linear on an open connected set $\Omega\subset\mathbb{R}^n$ if and only if $\Omega\subset\textrm{}^c\mathcal{A}_p.$ The injective function $\nu: \{E\} \rightarrow \mathbb{Z}_n\cap\mathcal{N}_p$ is induced by the order of the component. 

\begin{definition} \rm (See~\cite{Forsberg-Passare-Tsikh:2000}).
{\it The order of a connected component} $E\subset\mathrm{}^c\mathcal{A}_p$ is defined to be the vector $\nu\in\mathbb{Z}^n$ with the components $$\nu_j=\frac{1}{(2\pi i)^n}\int_{\mathrm{\scriptsize Log}^{-1}(u)}\frac{z_j\partial_j p(z)}{p(z)}\frac{dz_1\wedge\ldots\wedge dz_n}{z_1\ldots z_n},\quad j=1,\ldots,n,$$ where $u\in E.$ Values of~$\nu_j$ do not depend on the choice of~$u,$ since the homology class of the cycle~$\mathrm{Log}^{-1}(u)$ is the same for all $u\in E.$
\end{definition}

It is also possible to construct the compactified version of a polynomial amoeba by using the moment map~\cite{Mikhalkin:2000} $\mu: \left(\mathbb{C^*}\right)^n \rightarrow \mathcal{N}_p:$ $$\mu: (z_1, \ldots, z_n) \longmapsto \frac{\sum_{\alpha\in A}|z^\alpha|\cdot\alpha}{\sum_{\alpha\in A}|z^\alpha|}.$$

\begin{definition} \rm
{\it The compactified amoeba} $\bar{\mathcal{A}}_{p}$ of a polynomial $p(z)$ is the closure of the image of its zero locus under the map $\mu(z).$
\end{definition}

In addition to polynomial amoebas in this survey we also consider objects closely related to them such as coamoebas, contours, and spines of amoebas, since computing these objects can often be helpful in the understanding of the properties of amoebas.

\begin{definition} \rm
{\it The coamoeba} $co\mathcal{A}_{p}$ of a polynomial $p(z)$ is the image of its zero locus under the map $\mathrm{Arg}:  \left(\mathbb{C^*}\right)^n \rightarrow \mathbb{R}^n:$ $$\mathrm{Arg}: (z_1, \ldots, z_n) \longmapsto (\mathrm{arg}(z_1), \ldots, \mathrm{arg}(z_n)).$$ 
\label{def:coamoeba}
\end{definition}
Coamoebas were introduced by M. Passare in a talk in 2004 and they have applications, for example, in the field of mathematical physics, see Section~\ref{sec:history_applications} for more detail. An analogue of the Ronkin function for coamoebas has been introduced in~\cite{Johansson-Samuelsson:2017}.

The $\mathrm{Log}$ and $\mathrm{Arg}$ maps generate the following commutative diagram~\cite{Sadykov-Tsikh-book}, where $P$ is the projection map onto the $n$-dimensional real torus:
 
\begin{figure}[ht!]\label{avb-intro}
$$
\xymatrix{
                    &\mathbb{(C^\ast)}^n\ar[dl]_{\text{\scriptsize Log}}  \ar[dr]^{\mathrm{Arg}}  \ar@/^/@{.>}[drr] & \\
            {\mathbb{R}^n} &  {\mathbb{C}^n}\ar[u]_{\mathrm{Exp}} \ar[l]^{\mathrm{Re}} \ar[r]_{\mathrm{Im}}  &
               {\mathbb{R}^n} \ar[r]_{\!\!\!\!\!\!\!\!\!\!\!\!\!\!\!P} & {(\mathbb{R}/2\pi\mathbb{Z})^n}.}
$$
\end{figure}

%
%


\begin{definition} \rm
{\it The contour} of the amoeba~$\mathcal{A}_{p}$ is the set~$\mathcal{C}_{p}$ of critical points of the logarithmic map $\mathrm{Log}$ restricted to the zero locus of the polynomial $p(x).$ 
\end{definition}

The contour of an amoeba is the closed real-analytic hypersurface in~$\mathbb{R}^n, $ the boundary $\partial\mathcal{A}_p$ is a subset of the contour~$\mathcal{C}_p$ but is in general different from it. Results on the maximal number of intersection points of a line with the contour of a hypersurface amoeba are given in~\cite{Lang-Shapiro-Shustin:2021}.

Let $\mathcal{A}'\subset\mathbb{R}^n\cap\mathcal{N}_p$ be a set of vectors $\alpha$ such that $\textrm{}^c\mathcal{A}_p$ contains components~$E_\alpha$ of the order~$\alpha$ and $$a_\alpha=\frac{1}{(2\pi i)^n}\int_{\mathrm{\scriptsize Log}^{-1}(u)}\mathrm{ln}|\frac{p(z)}{z^\alpha}|\frac{dz_1\wedge\ldots\wedge dz_n}{z_1\ldots z_n},\quad\forall\alpha\in\mathcal{A}',\quad u\in E_\alpha.$$

Consider the function $S_p:\mathbb{R}^n\rightarrow\mathbb{R}:$

$$S_p(x)=\max\limits_{\alpha\in\mathcal{A}'}\{\langle\alpha,x\rangle+a_\alpha\}$$

\begin{definition} \rm
The set $\left\{x\in\mathbb{R}^n,\textrm{ where }S_p(x)\textrm{ is nonsmooth}\right\}$ is called {\it the spine of}~$\mathcal{A}_p.$
\label{def:spine}
\end{definition}

It is proved in~\cite{Passare-Rullgard:2000} that the spine of the amoeba~$\mathcal{A}_p$ is its deformation retract.

\begin{example} \rm
Consider the polynomial $p_1(z_1,z_2)=5 z_1 + 15 z_1^2 + 8 z_1^3 z_2 + 10 z_2^2 + 10 z_1^3 z_2^2 + 8 z_2^3 + 15 z_1 z_2^4 + 5 z_1^2 z_2^4 + 50 z_1 z_2^3 + 50 z_1^2 z_2.$ The Newton polytope, the amoeba, its spine, and the compactified amoeba of~$p_1$ are shown in Figure~\ref{fig:p1_NewtonPolytope}.

\begin{figure}[h!]
\hskip0.5cm
\begin{minipage}{4.5cm}
\includegraphics[width=4.5cm]{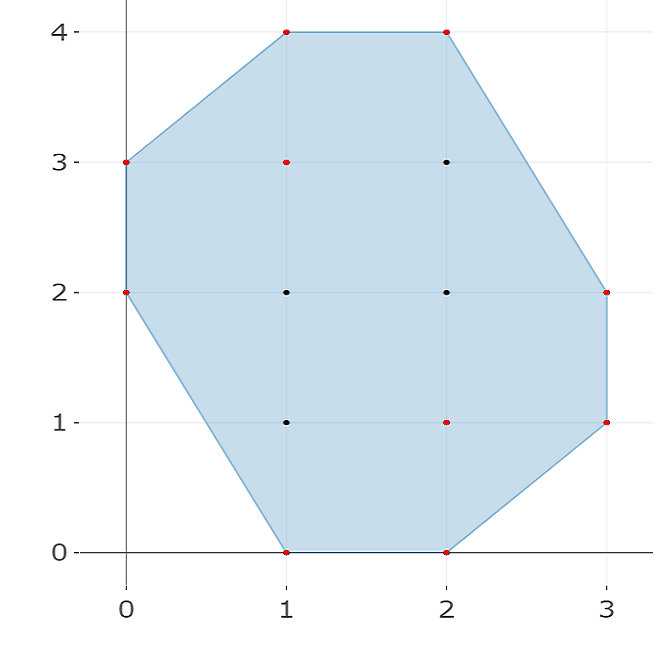}
\centering $(a)$
\end{minipage}
\hskip1cm
\begin{minipage}{5cm}
\includegraphics[width=5cm]{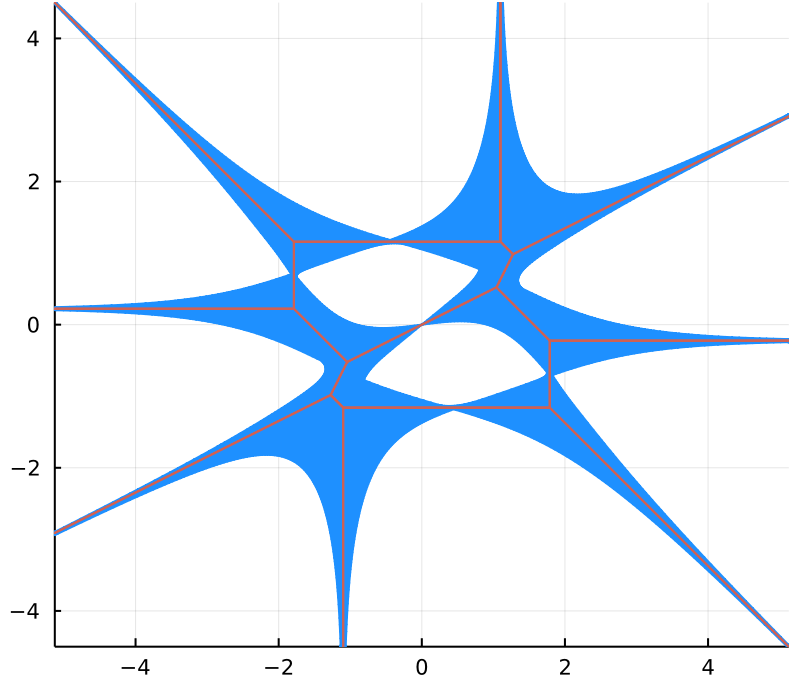}
\centering $(b)$
\end{minipage}

\begin{minipage}{5.5cm}
\includegraphics[width=5.5cm]{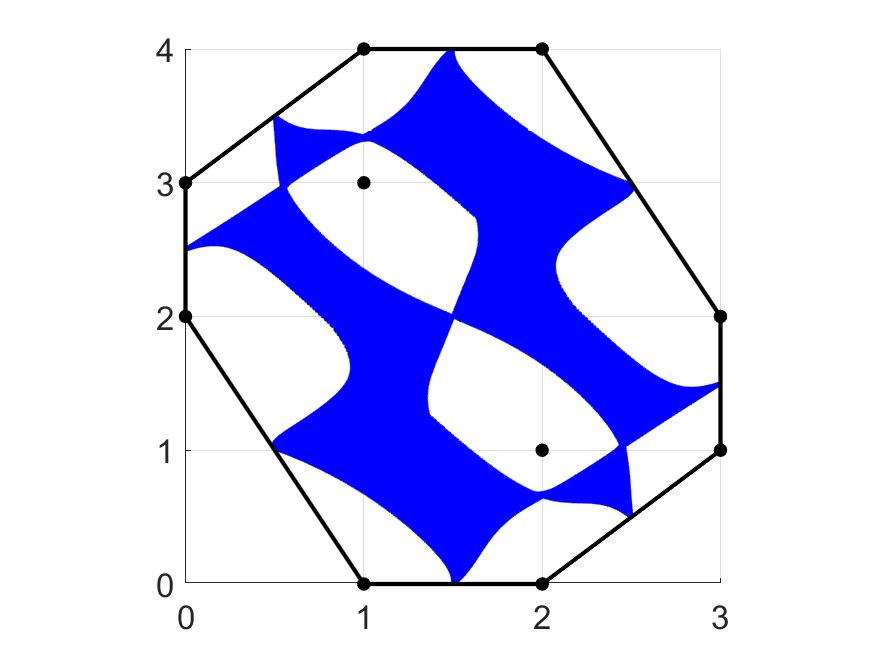}
\centering $(c)$
\end{minipage}
\begin{minipage}{5cm}
\includegraphics[width=5cm]{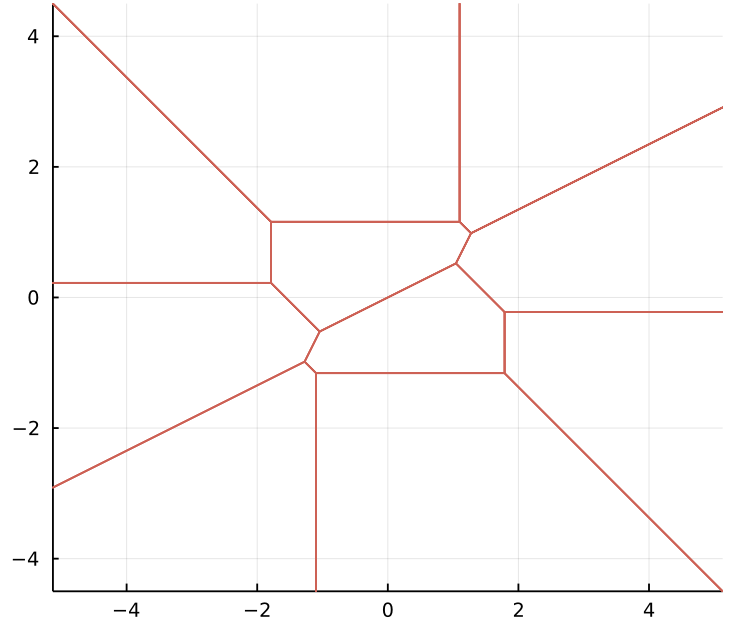}
\centering $(d)$
\end{minipage}
\caption{The Newton polytope~(a), the amoeba~(b), the compactified amoeba~(c), and the spine of the amoeba~(d) of the polynomial~$p_1$}
\label{fig:p1_NewtonPolytope}
\end{figure}
\end{example}

\subsection*{The Passare conjecture}

An important question whether one should try to prove the Passare conjecture~\ref{conj:Passare} or to refute it is opened at the moment. M. Passare stated that ``it would seem very plausible that the number of complement components is minimal for maximally sparse polynomials with at most $n + 2$ terms''~\cite{PassareConjectureFormulation:2009}. Numerous computational experiments imply that in lower dimensions (in particular, for $n=2$) maximally sparse polynomials have solid amoebas.

\begin{example} \rm
Consider the polynomial~$\tilde{p}_1(z_1, z_2)$ obtained by dropping from~$p_1(z_1, z_2)$ the monomials corresponding to the inner lattice points of its Newton polytope: $\tilde{p}_1(z_1, z_2) = 5 z_1 + 15 z_1^2 + 8 z_1^3 z_2 + 10 z_2^2 + 10 z_1^3 z_2^2 + 8 z_2^3 + 15 z_1 z_2^4 + 5 z_1^2 z_2^4.$ The amoeba of $\tilde{p}_1$ is shown in Figure~\ref{fig:p1_solid_Amoeba} and it is solid.

\begin{figure}[h!]
\centering
\begin{minipage}{5cm}
\includegraphics[width=5cm]{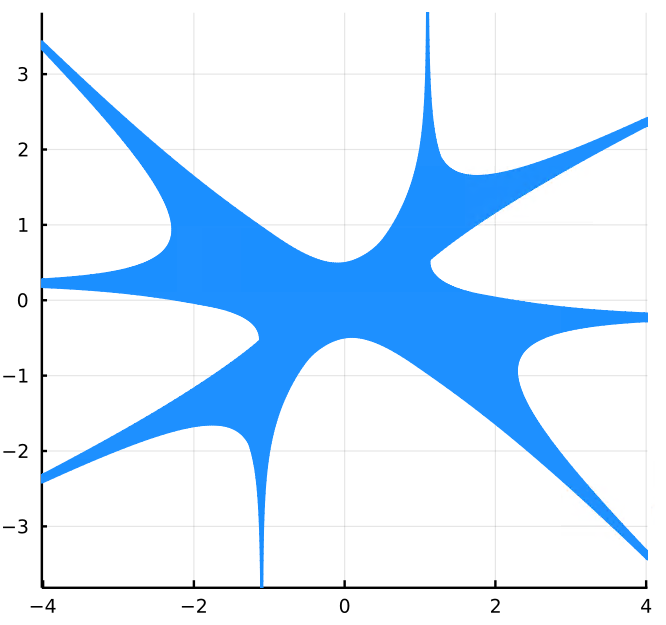}
\end{minipage}

\caption{The amoeba of the polynomial $\tilde{p}_1$ is solid}
\label{fig:p1_solid_Amoeba}
\end{figure}
\end{example}

There are results on solid amoebas for some classes of sparse polynomials. In~\cite{Iliman-deWolff:2016} the authors consider the class of polynomials of several real variables whose Newton polytopes are simplices and their supports are given by all the vertices of the Newton polytopes and a single additional interior lattice point. It is proved that under some conditions amoebas of such polynomials are solid. Another family of polynomials with solid amoebas, as proven in~\cite{Passare-Sadykov-Tsikh:2005}, are principal $\mathcal{A}$-determinants introduced in~\cite{GKZ:1994}, that also implies solidness of amoebas of discriminants of univariate polynomials. 

In the context of the Passare conjecture the theory of fewnomials developed by A.G. Khovanskii~\cite{Khovanskii:1991} should be mentioned. The whole ideology behind it is that the variety defined by a system of ``simple'' polynomials should have a ``simple'' topology. A maximally sparse polynomial is simple in the same way, that is, it contains the minimal possible number of monomials among all of the polynomials with a~given Newton polytope. 

The Passare conjecture raises the question of the dependence between the support of a polynomial and connected components in the complement of its amoeba. In particular, it could be easily proven if, similarly to the order of component map $\nu,$ there existed an injective map from $\{E\}$ to $A.$ The following example shows that this is not the case.

\begin{example} \rm
Consider the polynomial~$p_2(z_1,z_2)=z_1^2+50z_2^3+100iz_1z_2^3+100z_1^3z_2^3-50z_1^4z_2^3+z_1^2z_2^6.$ The Newton polytope and the amoeba of $p_2$ are given in Figure~\ref{fig:extra_component}. 
\begin{figure}[h!]
\hskip0.5cm
\begin{minipage}{4cm}
\includegraphics[width=4cm]{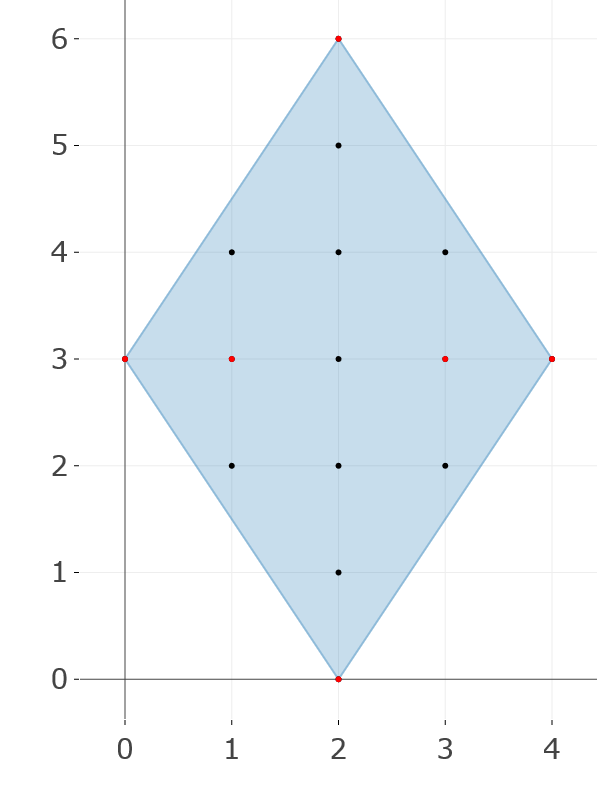}
\centering $(a)$
\end{minipage}
\hskip1cm
\begin{minipage}{6cm}
\includegraphics[width=6cm]{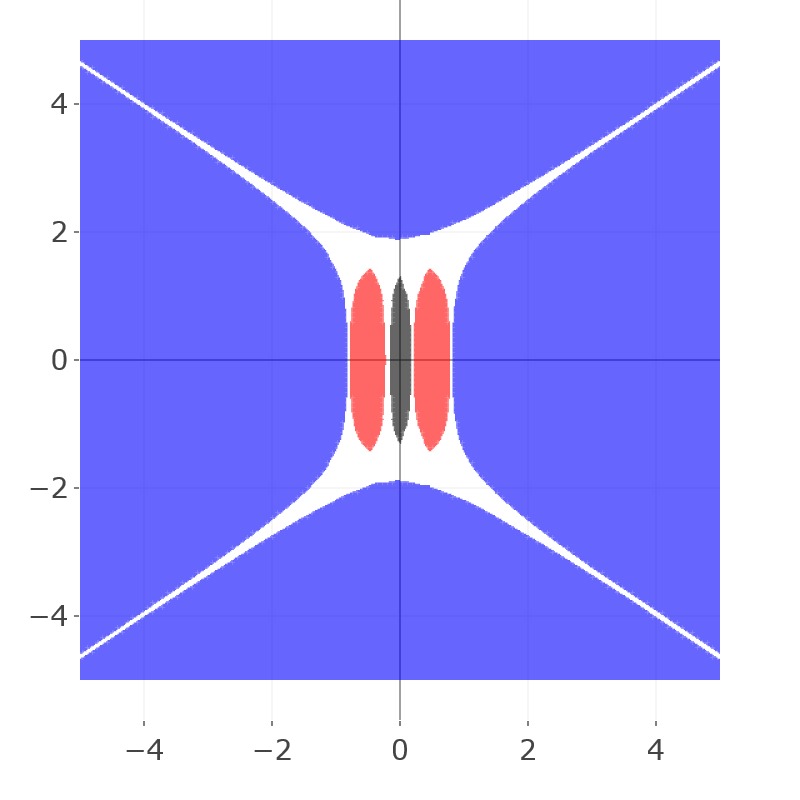}
\centering $(b)$
\end{minipage}
\caption{The Newton polytope (a) and connected components of the complement of the amoeba (b) of the polynomial~$p_2$}
\label{fig:extra_component}
\end{figure}
\label{ex:extra_component}
\end{example}

The unbounded components of the amoeba complement correspond to the vertices of the Newton polytope. The polynomial~$p_2$ includes another two monomials with the exponent vectors~$(1,3)$ and~$(3,3),$ but there are three bounded components in the amoeba complement. Thus there is the bounded component in the amoeba complement which does not correspond to any monomial.

It should be also noted that the dependence between the coefficients of the polynomial with a~given Newton polytope and the number of bounded components in its amoeba complement is quite intricate. The following example illustrates one of the problems arising from this fact.

\begin{example} \rm
\label{ex:monstera}
In Figure~\ref{fig:monstera} (a) a~leaf of a~tropical plant called monstera is shown, the picture is a~courtesy of T.M.~Sadykov. 
\begin{figure}[h!]
\hskip0.5cm
\begin{minipage}{5.5cm}
\includegraphics[width=5.5cm]{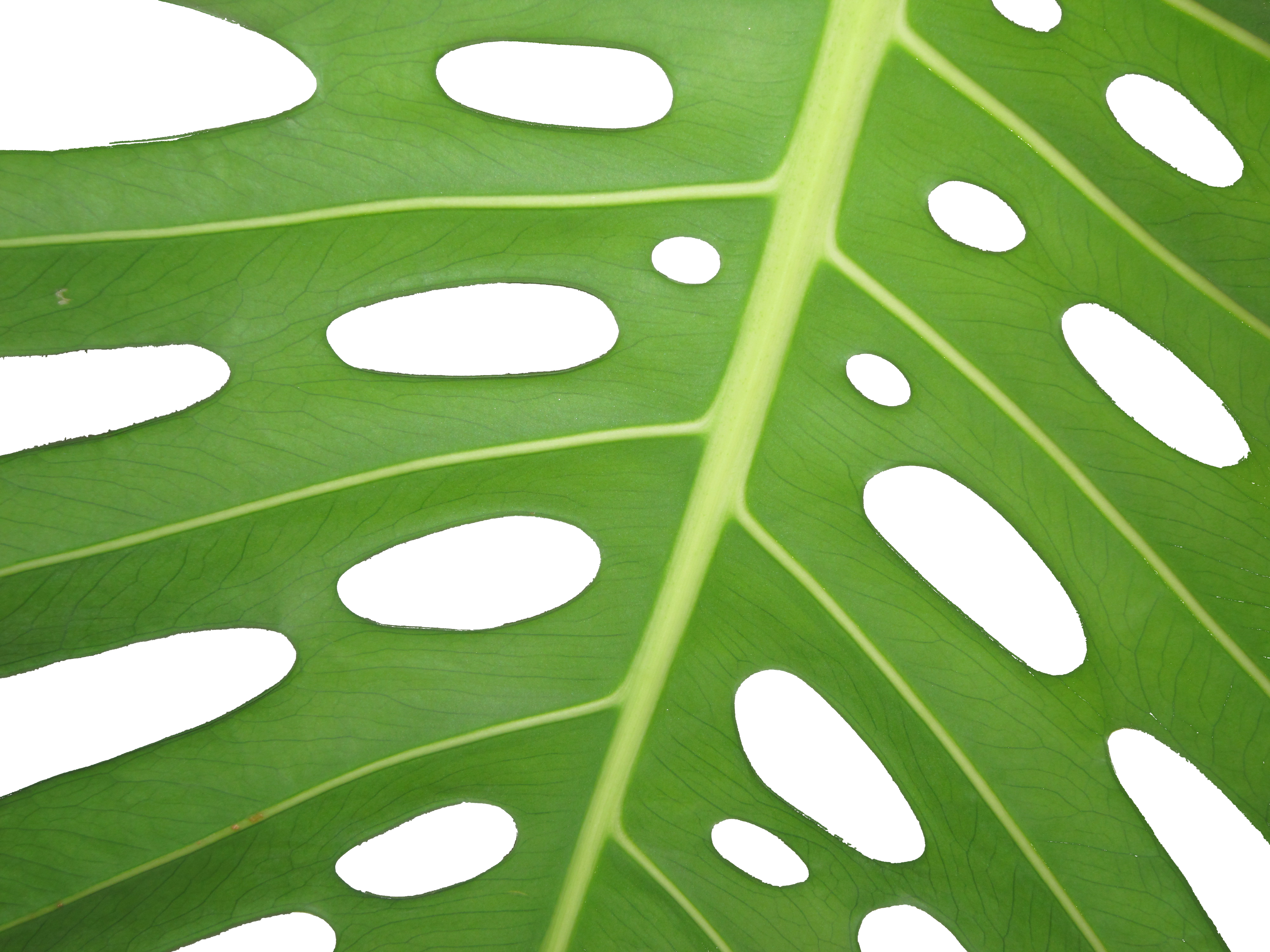}
\centering $(a)$
\end{minipage}
\hskip0.5cm
\begin{minipage}{7cm}
\includegraphics[width=7cm]{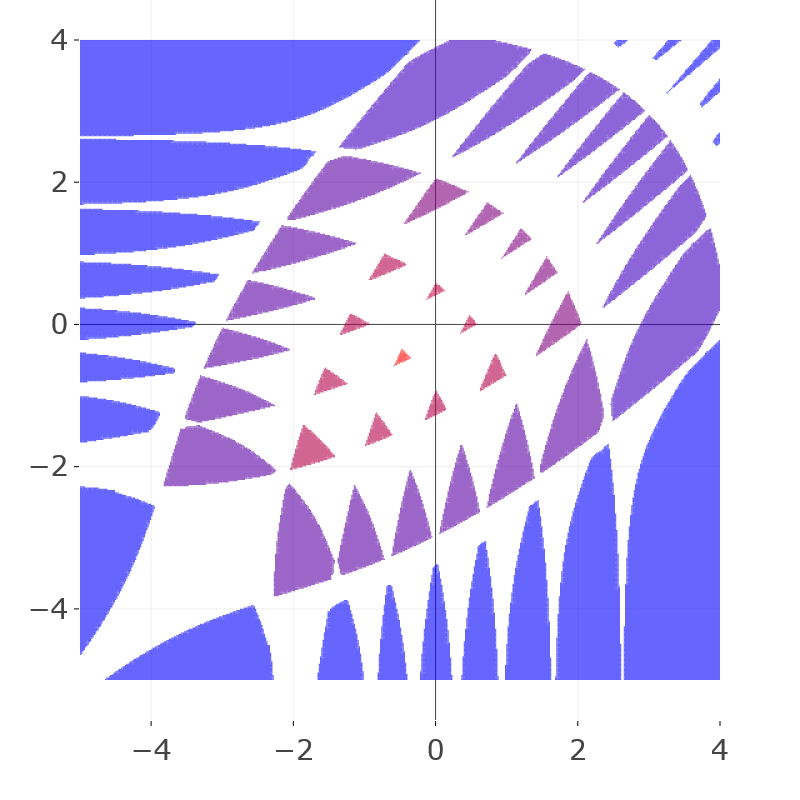}
\centering $(b)$
\end{minipage}

\caption{A monstera leaf (a) and connected components of the complement of its amoeba counterpart (b)}
\label{fig:monstera}
\end{figure}
Note the similarity in the structure between this leaf and the amoeba of a~bivariate polynomial. Holes in the monstera leaf are convex, just like the connected components of the polynomial amoeba complement, and parts of the leaf with cuts resemble the tentacles of an amoeba. Openings of this shape allow the plant to optimize the process of photosynthesis on the entire surface of the leaf and leaves that are in the shade often do not have such cuts. The pattern of veins on the leaf resembles the structure of the amoeba spine (see Definition~\ref{def:spine}).

In Figure~\ref{fig:monstera} (b) is depicted the closest amoeba counterpart we could find -- namely, it is the optimal amoeba of the polynomial generated by the means of the \texttt{OptimalPolynomialSimplex} algorithm from~\cite{Zhukov-Sadykov:2023}. The full form of the polynomial is omitted here due to its rather large size, one can obtain it by using the algorithm with the parameters $\{$dimension = 2, degree = 10, power = 2$\}$ and dropping the monomial corresponding to the exponent vector $(0,0)$ from the result.

The similarity of two objects is rather superficial and this leads to the problem of finding the explicit form of a polynomial given its amoeba. In particular cases, it can be solved by choosing the coefficients of a polynomial with a given Newton polytope, but in general this process is not so simple, since a change in a single coefficient can lead to the unpredictable appearance or disappearance of components in the amoeba complement. Other ``tropical'' objects related to the topic of polynomial amoebas are discussed in Section~\ref{sec:tropical_geometry}

\end{example}

\section{Historical reference and applications of polynomial amoebas}
\label{sec:history_applications}
The reason for the great interest in computing polynomial amoebas and studying their properties is the fact that amoebas characterize the zero loci of polynomials and thus, on the one hand, have been the object of study even before the actual concept of amoeba appeared, and on the other hand, they are widely used in various domains of science. This section discusses mathematical concepts related to polynomial amoebas, both historical and recent ones, as well as applications of polynomial amoebas.

\subsection{Amoebas and convergence of power series}

Consider a power series centered at the origin:
\begin{equation}
\sum\limits_{\alpha\in\mathbb{N}^n}c_\alpha z^\alpha, c_\alpha\in\mathbb{C}, z^\alpha=z_1^{\alpha_1}z_2^{\alpha_2}\ldots z_n^{\alpha_n}.
\label{eq:powerSeries}
\end{equation}

The domain of convergence of the series~(\ref{eq:powerSeries}) is a complete Reinhardt domain centered at the origin. For $n = 1$ it is stated in Abel’s lemma and attempts to find the original wording of this fact have led to the formulation by A.-L. Cauchy in 1831 \cite[vol. 8, p. 151]{Cauchy-OeuvresCompletes} (translated from the French original):


{\it $x$ designating a real or complex variable, a real or complex function of $x$ can be developed into an ordered convergent series according to the ascending powers of this variable, as long as the modulus of the variable retains a value lower than the smallest of those for which the function or its derivative
ceases to be finite or continuous.}

The proof of this statement in the case of multiple variables is given in several different sources, but it is not an easy task to detect the initial one.

If all of the singularities of the series~(\ref{eq:powerSeries}) belong to an algebraic hypersurface then the structure of its domain of convergence is closely related to the amoeba of this hypersurface. More precisely, for any polynomial $p(z)$ there exists a bijective correspondence between connected components of the
amoeba complement for an algebraic set $\{p(z) = 0\}$ and domains of convergence of the Laurent series expansions with denominator $p(z)$~\cite{GKZ:1994}.

\subsection{Hilbert's 16th problem}

One of the problems related to polynomial amoebas is the famous Hilbert’s 16th problem on the topology of algebraic curves and surfaces~\cite{Hilbert:1902}. Originally it was formulated as consisting of two parts: the first part considering the relative positions of the branches of real algebraic curves of degree $d$ and the second one on the upper bound for the number of limit cycles in two-dimensional polynomial vector fields of degree $d$ and their relative positions.

The first part of Hilbert’s 16th problem was followed by A. Harnack’s investigation on the number of separate connected components for algebraic curves of degree $d$ in $\mathbb{R}\mathbb{P}^2$. Harnack proved~\cite{Harnack:1876} that the number of components of such a curve does not exceed $\frac{(d-1)(d-2)}{2}+1$. If the maximal number of components is attained for a curve then it is called a {\it Harnack curve}. The components of a Harnack curve have the best possible topological configuration in some sense~\cite{Mikhalkin:2000}. The amoeba of a Harnack curve can be described analytically as $$\left\{(x,y)\in\mathbb{R}^2: \prod\limits_{\alpha,\beta=\pm 1} P_W(\alpha e^x,\beta e^y)\leq 0\right\},$$
where $P_W$ is the polynomial defining the curve~\cite{Feng-He-Kennaway-Vafa:2008}.

Amoebas of Harnack curves enjoy different optimal properties: for example, the amoeba of Harnack curve with genus $g$ has exactly $g$ compact components in its complement~\cite{Kenyon-Okounkov-Sheffield:2006}, a curve is Harnack if and only if its amoeba has the maximal possible area for a given Newton polygon~\cite{Mikhalkin-Rullgard:2001}, a Harnack curve~$C$ always possesses the map from $C(\mathbb{C})$ to $\mathcal{A}_C$ such that it is 2-to-1 except for, maybe, a finite number of real nodes where it is 1-to-1.  

The special kind of variational principle for the amoeba ``holes'' (that is, bounded components of the amoeba complement) is formulated in~\cite{Kenyon-Okounkov:2006}, similar to the variational principle for conformal maps and it is proved in this work that the areas of the amoeba holes and the distances between the amoeba tentacles provide global coordinates on the moduli space of Harnack curves.

An explicit integral formula is presented in~\cite{Passare:2016} for the amoeba-to-coamoeba mapping in the case of polynomials defining Harnack curves as well as the exact description of the coamoebas of such polynomials.


%
%

\subsection{Modern amoebas}

There are numerous recent papers on different topics related to polynomial amoebas. Amoebas can be defined for different varieties including spherical tropical varieties~\cite{Kaveh-Manon:2019}, non-hypersurface ones~\cite{Juhnke-Kubitzke-deWolff:2016}, and subvarieties of non-commutative Lie groups~\cite{Mikhalkin-Shkolnikov:2022}. Zero-dimensional case is considered in~\cite{Nisse:2016}. Results on the dimension of amoebas of varieties are presented in~\cite{Draisma-Rau-Yuen:2020,Mikhalkin:2017}. The special kind of the polyhedral complex which is the subset of the Newton polytope of a polynomial and enjoys the key topological and combinatorial properties of the polynomial amoeba is considered in~\cite{Nisse-Sadykov:2019}.

Some of the recent articles reformulate classical problems in terms of the polynomial amoeba theory. For example, in~\cite{Theobald-deWolff:2016} this approach is applied to the problem of defining the behaviour of univariate trinomial roots. In~\cite{Lyapin:2009} polynomial amoebas are used for investigating the properties of Riordan's arrays which arise as solutions of Cauchy problem for difference equations.

Polynomial amoebas enjoy some optimal properties if they are defined by hypergeometric functions. It is proved in~\cite{Bogdanov-Sadykov:2020} that under certain nondegeneracy conditions amoebas of hypergeometric polynomials are optimal. In~\cite{Passare-Sadykov-Tsikh:2005} it is shown that the singular hypersurface of any nonconfluent hypergeometric function has a solid amoeba. The description of the convergence domain of a hypergeometric series in terms of the amoeba complement is given in~\cite{Passare-Tsikh:2002}, the more recent results on the topic are presented in~\cite{Cherepanskiy-Tsikh:2020,Nilsson-Passare-Tsikh:2019}.

\subsection{Tropical geometry}
\label{sec:tropical_geometry}
A lot of natural connections to polynomial amoebas lie in the field of tropical geometry. Originally the word ``tropical'' refers to one of the pioneers of this theory -- Brazilian computer scientist and logician I.~Simon~\cite{Yger:2012}. In this survey, we will not dive deep into the specifics of this area, since there are many researches, both classical and recent, that describe in great detail the relationship between polynomial amoebas and tropical geometry (see, for example, \cite{Kim-Nisse:2021,Lang:2020,Mikhalkin:2004,Mikhalkin:2006,deWolff:2017,Yger:2012}). However, since some of the algorithms for computing amoebas involve the use of tropical geometry concepts, some introduction is still required.

The basic object of tropical geometry is the {\it tropical semiring} $(\mathbb{R}\cup \{\infty\}, \oplus, \odot),$ which is basically a set of real numbers with infinity and operations of addition and multiplication being defined as follows: $x \oplus y = \textrm{max} (x,y), \quad x\odot y = x + y.$ Various frames for developing calculus over this structure are described in~\cite{Yger:2012}.

For the variables $x = (x_1,x_2,\ldots,x_n), x_i\in (\mathbb{R}\cup \{\infty\}, \oplus, \odot), i=1,\ldots,n$ we define {\it tropical polynomial}:

$$\mathfrak{p}(x_1,\ldots,x_n) = \bigoplus\limits_{k=1}^d a_{j_k} \odot x^{\odot j_k} = \textrm{max}\{a_{j_1}+\langle x, j_1\rangle, \ldots, a_{j_d}+\langle x, j_d\rangle\}.$$
A tropical polynomial in~$n$ variables is a piecewise-linear concave function on $\mathbb{R}^n$ with integer coefficients. {\it Tropical hypersurface} $\mathcal{H}(\mathfrak{p})$ is the set of points~$x\in \mathbb{R}^n,$ where $\mathfrak{p}(x)$ is not a linear function.  An example of the tropical hypersurface is the spine of an amoeba, since in terms of the tropical semiring, the function $S_p(x)$ from Definition~\ref{def:spine} is a tropical polynomial: $S_p(x) = \bigoplus\limits_{\alpha\in\mathcal{A}'} \odot x^{\odot \alpha}.$

Another connection of polynomial amoebas and tropical geometry is the tropicalization of an algebraic curve~\cite{Jonsson:2016}. An implementation of numerical methods for constructing the tropicalization of a complex curve is presented in~\cite{Jensen-Leykin-Yu:2016}.

There are algorithms for polynomial amoeba computation based on the notion of Archimedean amoebas~\cite{Avendano-etAl:2018, Forsgard:2020} (see Section~\ref{sec:algorithms} for the details). The topic of non-Archimedean amoebas is discussed in~\cite{Einsiedler-Kapranov-Lind:2006}.

\subsection{Applications in Physics}

There are two-dimensional combinatorial systems defined in the domain of mathematical physics, called dimer models~\cite{Kenyon:2008}. Dimer models
have vast applications within mirror symmetry and string theory. In some cases, dimer model graphs for curves in $\mathbb{C}^3$ agree with amoebas for these curves~\cite{Feng-He-Kennaway-Vafa:2008}. In~\cite{Kenyon-Okounkov-Sheffield:2006} the spectral curve of the Kasteleyn operator of the graph is considered. The amoeba of the spectral curve represents the phase diagram of the dimer model. One of the physical interpretations for dimer models are crystal surfaces at equilibrium, and such models with crystal facets being in bijection with the components of the complement of the amoeba are presented in~\cite{Kenyon-Okounkov-Sheffield:2006}. Also a lot of connections are established in this work between properties of amoebas of Harnack curves and phases of the dimer models (frozen, liquid, and gaseous). 

Not only amoebas but also coamoebas (see Definition~\ref{def:coamoeba}) have their applications in dimers theory. The relationship between dimer models on the real two-torus and coamoebas of curves in $(\mathbb{C}^\times)^2$ is described in~\cite{Forsgard:2019}. It is shown in this work that the dimer model obtained from the shell of the coamoeba is a deformation retract of the closed coamoeba if and only if the number of connected components of the complement of the closed coamoeba is maximal.

Another object related to polynomial amoebas is $(p, q)$-web model in field theory introduced in~\cite{Aharony-Hanany-Kol:1998}. The spine of the amoeba of a bivariate polynomial $P$ is the $(p, q)$-web associated to the toric diagram which is the Newton polygon of $P$~\cite{Feng-He-Kennaway-Vafa:2008}. In~\cite{Bao-He-Zahabi:2022} authors adapt the Mahler measure from number theory to toric quiver gauge theories and the Mahler flow they introduce there can be geometrically interpreted in terms of polynomial amoebas and their holes.

A solution to the system of fundamental thermodynamic relations in statistical thermodynamics leads to the notions of the polynomial amoeba and its contour~\cite{Passare-Pochekutov-Tsikh:2013}. It is shown in~\cite{Zahabi:2021} by using the quiver gauge theory that thermodynamic observables such as free energy, entropy, and growth rate are explicitly derived from the limit shape of the crystal model, the boundary of the amoeba, and its Harnack curve characterization. In~\cite{Konopelchenko-Angelelli:2018}, while studying the partition functions in different branches of physics, the authors introduce the notion of statistical amoebas and describe their relation with polynomial amoebas.

\section{Computation of polynomial amoebas}

The problem of giving a complete geometric or combinatorial description for the amoeba of a polynomial has a significant computational complexity, especially for the higher dimensions. This section contains a review on the main problems of the amoeba computation and the existing algorithms for their solution.

The implementation of algorithms can depend on different factors such as the programming languages or frameworks used, in particular, this applies to the data types. In this survey, technical details of the implementation are omitted where it is possible.

\subsection{Problems and algorithms}

\subsubsection{Membership problem}

One of the basic problems of amoeba computation is whether a given point belongs to the amoeba or, equivalently, if it belongs to a component of the amoeba complement (the membership problem). This problem is very natural, so numerous papers considering the computation of amoebas address it. One way to provide a solution for the membership problem is finding a certificate~$C$ such that if $C(|z|)$ is true then $\mathrm{Log} (z)\notin \mathcal{A}_p$ for a polynomial~$p$ (see, for example,~\cite{Forsgard-Matusevich-Mehlhop-deWolff:2018}). It is stated in~\cite[Corollary~2.7]{Theobald:2002} that for a fixed dimension~$n$ this problem can be solved in polynomial time. In~\cite[Theorem 4.2]{Avendano-etAl:2018} it is shown that in general this problem is {\bf PSPACE}, that is, it can be solved in polynomial time by a~parallel algorithm, provided one allows exponentially many processors.

One of the amoeba properties simplifying the computations is the lopsidedness (see~\cite{Forsgard:2021}). In~\cite{Purbhoo:2008} the lopsidedness criterion is presented, which provides an inequality-based certificate for non-containment of a point in the amoeba. Based on this concept, a converging sequence of approximations for the amoeba can be devised. These approximations use {\it cyclic resultants}, a fast method of computing such resultants is provided in~\cite{Forsgard-Matusevich-Mehlhop-deWolff:2018}. Another possible base for approximations like these is the real Nullstellensatz (see~\cite{Theobald-deWolff:2015}).

To solve the membership problem more efficiently, some authors propose formulations for it, that are less restrictive in some way. For example, in~\cite{Timme-Master} the soft membership problem is given it the following formulation: to determine whether a given point belongs to the amoeba {\it with high confidence}, which means that the chosen criterion should only provide the correct solution with some controllable high probability. This new problem is then reduced to the solution of a system of polynomial equations by means of {\it realification} of the initial polynomial (see~\cite{Theobald-deWolff:2015}). 

\subsubsection{Depiction of amoebas}
\label{sec:algorithms}
In two- and three-dimensional cases one of the simplest ways of describing the structure of an amoeba is depicting it. For $n > 3$ it is possible to visualize an amoeba by depicting its three-dimensional sections.

Since the tentacles of an amoeba $\mathcal{A}_p$ always stretch to infinity, only a part of the amoeba is usually depicted by choosing a domain $\Omega \subset \mathbb{R}^n$ and depicting the intersection $\mathcal{A}_p \cap \Omega.$ This fact raises the computational problem of choosing the domain $\Omega$ such that the intersection inherits the essential characteristics of the amoeba like the number of connected components in its complement. A~heuristic algorithm for computing $\Omega$ is presented in~\cite{Timme-Master}. 



%

The basic algorithm for depicting the amoeba of a polynomial includes choosing a grid on~$\Omega$ and reducing the problem to finding the roots of a univariate polynomial. For example, in the two-dimensional case for the polynomial~$p(z_1,z_2),$ the domain $\Omega=[x_{\mathrm{min}}, x_{\mathrm{max}}]\times[y_{\mathrm{min}}, y_{\mathrm{max}}],$ the grid~$\{(x_j, \theta_k) | (j,k)\in S \subset \mathbb{N}^2\}$ is chosen on the set~$[x_{\mathrm{min}}, x_{\mathrm{max}}]\times[0,2\pi)$ and for each grid point $(x_j, \theta_k)$ we compute the roots of the polynomial~$p(e^{x_j+i\theta_k},z_2).$ For each of the roots~$z_2^{(m)}$ by the definition of the polynomial amoeba $(x_j, \mathrm{ln}|z_2^{(m)}|) \in \mathcal{A}_p.$ Such algorithms in different variations are presented in~\cite{Bogdanov-Kytmanov-Sadykov:2016,Forsberg-PHD,Johansson-PHD,Leksell-Komorowski-Bachelor,Nilsson-PHD}. We refer to these algorithms as ``naive'' ones in what follows. Naive algorithms have some major drawbacks. First of them is that not all of the amoeba points obtained with such an algoritm belong to~$\Omega,$ and due to this fact some of the performed computations become unnecessary. The good quality of the picture can require a~grid with a~very small step because of the low control on the density of obtained points of the amoeba. Another feature of these algorithms is that each of them includes as a subroutine an algorithm for computing the roots of a univariate polynomial with complex coefficients. There are numerous algorithms for solving the root-finding problem, some of them are listed in~\cite{Leksell-Komorowski-Bachelor}.  

In~\cite{Timme-Master} another approach to depicting amoebas is presented, based on the mapping of the domain~$\Omega$ into a set of pixels of an output device and executing  the membership test for corresponding points to determine whether a pixel should be depicted as a part of the amoeba or not. In what follows we refer to this algorithm as ``simple'' one as its author does. 

Some of the improvements to the computation of amoebas include construction of the grid based on Archimedean tropical hypersurfaces~\cite{Avendano-etAl:2018} which allow to execute membership tests only for the points lying in the close neighborhood of the amoeba.  Another way to avoid testing the points which  do not belong to the amoeba is by using greedy algorithms~\cite{Timme-Master}, that is, testing only points in some neighborhood of already tested ones. 

One more approach, also presented in~\cite{Timme-Master}, is based on the approximation of the amoeba for a bivariate polynomial with a set of polygons. The procedure starts with the construction of the amoeba spine and a set of triangles such that each of them has one side coinciding with one of the spine segments. Then new triangles are added iteratively so that their union covers the area of the amoeba as much as possible. The membership test is used to determine the size of new added triangles.

In~\cite{Zhukov-Sadykov:2023} authors compute the connected components of the amoeba complement by using their orders. Since order values are the same for all points in the same component and differ for the points from the different components, the order of the component is a good classifier for the points of the amoeba complement. For the points of the amoeba the integral for computing the order does not always converge but the jump of the order function itself can be a criterion for containment of a point in the amoeba. The construction process itself is similar to the polygonal algorithm from~\cite{Timme-Master}, but uses rectangular parallelepipeds instead of triangles and thus can easily be reproduced for polynomials with an arbitrary number of variables. Hereinafter we refer to this algorithm as ``dichotomous'' since it includes iterative division of parallelepipeds into $2^n$ parts to increase the precision of computations.

It must be noted that there are common techniques that allow improving the performance of algorithms for amoeba computation. One of the examples is the use of distributed computing, the parallel algorithms for computing amoebas are given in~\cite{Leksell-Komorowski-Bachelor, Timme-Master}.

In~\cite{Bao-He-Hirst:2023} computational problems for polynomial amoebas are solved by means of machine learning. Authors use artificial neural networks to determine the genus of an amoeba and to solve the membership problem. Some of algorithms there demonstrate impressive results, for example, one of the models for classification of amoebas based on their genus achieves the prediction accuracy around 0.95. By the accuracy here we denote the metric for evaluating classification models calculated as a~ratio of the number of correct predictions to the total number of predictions. Neural network models seem to be a suitable tool for solving some problems related to amoeba polynomials, for example, the problem of finding a polynomial given its amoeba from Example~\ref{ex:monstera} looks like an image classification problem that can be solved using neural networks.

\subsubsection{Problem of small components} 
\label{sec:small_components}
This problem is not explicitly formulated in the literature, however, several works on the visualization of polynomial amoebas address the problem of bounded complement components having a small size (``small'' or ``narrow'' components). 

To introduce a strict definition of the size of a bounded component of the amoeba complement, one can define the diameter of the component as the maximal distance between its points. In~\cite{Zhukov-Sadykov:2023} it is stated that for the dichotomous algorithm of the amoeba visualization detection of the components of the complement with the diameter less than some small~$\varepsilon$ is not guaranteed. In~\cite{Timme-Master} the algorithm for calculating the spine of an amoeba uses as one of the input parameters a number~$\delta$ that the author calls a lower bound of the size of any component of the complement and the runtime of this algorithm actually depends on~$\delta$.

The problem of small components can be crucial for finding a counterexample to the Passare conjecture and it is an open question whether it is possible to develop an algorithm that finds the smallest diameter of the bounded component of the amoeba complement in polynomial time with respect to~$n.$

\subsubsection{Algorithm comparison}

It should be noted that the computational complexity of the presented algorithms depends on the complexity of their subroutines, in particular, on the solution of the membership problem within their framework. Therefore, to simplify the presentation, the conclusions of the authors of the algorithms and the results of computer experiments are taken as a~basis for comparison.

Naive algorithms are best suited for understanding the general structure of an amoeba. Because of their simplicity compared to other algorithms, they allow one to compute amoebas relatively fast and accurate. However, if the amoeba of the input polynomial is sufficiently complex, for example, if its complement contains bounded components of a small size, such algorithms perform poorly.

The simple algorithm allows much better control of the computed amoeba points compared to the naive ones and thus it performs better on small complement components. The computation time for it is comparable to the results for naive algorithms, including the computations of optimal amoebas. Using of the grid based on Archimedean tropical hypersurfaces improves the computation time but not substantially~\cite[Section 1.4]{Timme-Master}.

A significant increase in the computation speed is given by using the greedy algorithm. Another important property of this algorithm is that an increasing number of the grid points (that is, computation with a higher precision) does not lead to the proportional increase of computation time~\cite[Remark 1.5.6]{Timme-Master}.

The polygonal approximation provides the best picture quality of all presented methods (see Figure~\ref{fig:p2_software_b} (b)). However, the frequent appearance of image artifacts is associated with this algorithm (see Example~\ref{ex:artifacts}), though it can be due to implementation peculiarities.

The dichotomous algorithm is comparable in the termination speed to the polygonal algorithm and it is also capable of performing computations with an arbitrarily high precision (though the increase in precision significantly increases the computation time). For detecting the small components with this algorithm, one can choose the domain $\Omega$ with an area small enough.

\subsection{Software comparison}

Let us consider three examples of the software tools for depicting polynomial amoebas and compare their functionality. All of these solutions are freeware and available online. 

{\bf [a]} The script by Dmitry Bogdanov at \url{http://dvbogdanov.ru/amoeba}. Very simple, both in usage and its functionality, it allows one to depict amoebas and compactified amoebas. The script itself only generates the MatLab code for depicting the amoebas with given parameters. Amoebas with the most complex structure amongst depicted using this software are presented at~\url{www.researchgate.net/publication/ 338341129_Giant_amoeba_zoo}.

{\bf [b]} The PolynomialAmoebas package written by Sascha Timme in the Julia programming language at \url{https://github.com/saschatimme/PolynomialAmoebas.jl}. It has very broad functionality, including not only tools for depicting amoebas and not only in two dimensions. It also allows one to depict spines and contours of amoebas, coamoebas, and amoebas in three dimensions. Computation of amoebas is possible by means of several algorithms.

{\bf [c]} The project by Timur Sadykov and Timur Zhukov at \url{http://amoebas.ru/index.html}. It includes the visualization of amoebas in two and three dimensions, three-dimensional slices of four-dimensional amoebas, and the evolution of the amoeba due to a change of the coefficients in the polynomial. 

\subsubsection{Functionality and algorithms}

For all of the presented software tools the basic functionality is depicting polynomial amoebas, but they differ significantly in terms of the implemented algorithms.

The software [a] implements naive algorithm for computing amoebas and its version for computing compactified amoebas. The drawbacks of the software are a lack of an algorithm that chooses a domain for depicting the amoeba and a necessary Matlab installation for performing computations.

The software [b] includes implementations of simple and greedy algorithms as well as algorithms based on Archimedean tropical hypersurfaces and polygonal approximation. The domain for depicting the amoeba can be automatically chosen by using the heuristic algorithm based on the computing of the Archimedean tropical hypersurface for the amoeba. Versions of simple and greedy algorithms are used for depicting three-dimensional amoebas. Computation of coamoebas is possible by the means of simple, greedy or coarse algorithms.

The software [c] implements the dichotomous algorithm based on the computing of the orders of the components of the amoeba complement. Due to the fact that the algorithm does not depend on the dimension of the space, in addition to two-dimensional amoebas, a depiction of three-dimensional amoebas and sections of four-dimensional ones is possible.

The summary of the functionality and implemented algorithms for computing amoebas in the software tools is presented in Table~\ref{tab:CompareFunctionality}. 

\begin{table}[h!]
\small
\caption{Functionality and implemented algorithms of the packages for computing polynomial amoebas}

\noindent\begin{tabular}{|p{6cm}|wc{2cm}|wc{2cm}|wc{2cm}|}
\hline
\multirow{2}{*}{Functionality} & \multicolumn{3}{c|}{Algorithms} \\
\cline{2-4}
& Software [a] & Software [b] & Software [c] \\
\hline
Newton polytopes & + & & + \\ 
\hline
Amoebas & naive & simple &  dichotomous \\
&&greedy&\\
&&Archimedean&\\
&&polygonal &\\
\hline
Compactified amoebas & + & & \\
\hline
Contours of amoebas & & + & \\
\hline 
Spines of amoebas & & + & \\
\hline
Coamoebas & &+& \\
\hline
3D amoebas & & simple & dichotomous \\
&&greedy&\\
\hline
Amoeba evolution & & & + \\
\hline
\end{tabular}
\vskip0.2cm
\label{tab:CompareFunctionality}
\end{table}


The important question for measuring the software performance is the choice of parameters. Depicting amoebas with the software [a] implies the choice of a grid by choosing numbers of modulus values~$r_x, r_y$ and a number of argument values~$\phi.$ Detailed information on the algorithm and its parameters is given in~\cite{Bogdanov-Kytmanov-Sadykov:2016}. For the lower values of these parameters, algorithm terminates faster, but the picture of amoeba dissipates into the set of isolated points, especially after zooming in it. For the tests, two sets of the parameters have been used -- the default ones (each number equals 100) and when each number equals 500, resulting in the smooth picture. For the simplicity in what follows these cases are denoted by ``100 values'' and ``500 values'' respectively.

Computation results with the software [b] depend on the choice of the algorithm and the grid. Of all the packages under consideration, this is the only one that provides not a graphical, but a~programming interface for accessing its functions.

For the software [c] the important parameter affecting the time of termination and the quality of resulting picture is the number of iterations the algorithm performs. The number of precision samples shows how many checks the algorithm performs on the order of points it finds (by computing the integral for different argument values). If this number is high, the algorithm terminates slower, but for lower values of this parameter some of the points of the amoeba can be misinterpreted for the points of the amoeba complement and vise versa.

\begin{example} \rm
Consider the polynomial~$p_3(z_1,z_2) = 5 z_1 + 15 z_1^2 + 240 z_1 z_2 + 400 z_1^2 z_2 + 8 z_1^3 z_2 + 10 z_2^2 + 900 z_1 z_2^2 + 900 z_1^2 z_2^2 + 10 z_1^3 z_2^2 + 8 z_2^3 + 400 z_1 z_2^3 + 240 z_1^2 z_2^3 + 15 z_1 z_2^4 + 5 z_1^2 z_2^4.$ This polynomial has been generated using the \texttt{OptimalPolyZonotope} algorithm given in~\cite{Zhukov-Sadykov:2023}. 
Let us compute the amoeba of~$p_3$ using the software [a], [b], and [c].

The result obtained with the software [a] for 100 values is shown in Figure~\ref{fig:p2_software_a} (a). It describes the structure of the amoeba but the quality of picture is not high enough. For 500 values the quality improves and there are no points in the picture that look isolated (see Figure~\ref{fig:p2_software_a} (b)). 

\begin{figure}[h!]
\hskip0.5cm
\begin{minipage}{5cm}
\includegraphics[width=5cm]{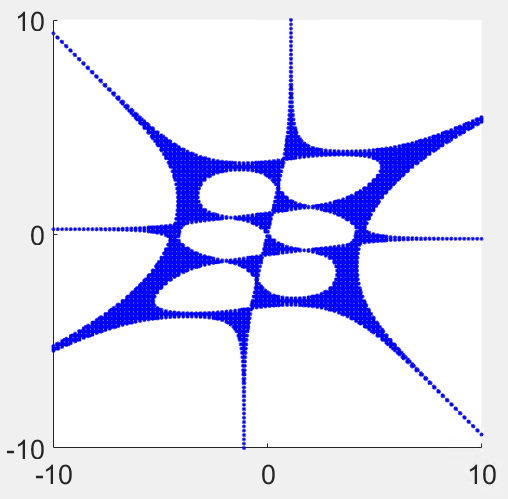}
\centering $(a)$
\end{minipage}
\hskip1cm
\begin{minipage}{5cm}
\includegraphics[width=5cm]{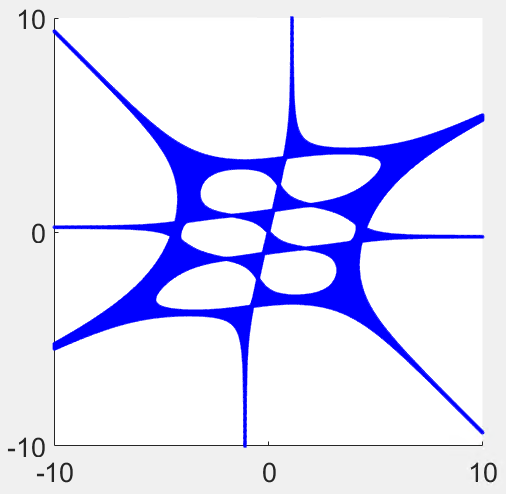}
\centering $(b)$
\end{minipage}

\caption{The amoeba of the polynomial~$p_3$ computed with the software [a] -- 100 values (a), 500 values (b)}
\label{fig:p2_software_a}
\end{figure}

The computation results for the software [b] are presented in Figure~\ref{fig:p2_software_b}. The fastest of the algorithms this package provides is the greedy one (see Figure~\ref{fig:p2_software_b} (a)), but the picture for the default algorithm (the polygonal approximation) has better quality (see Figure~\ref{fig:p2_software_b} (b)).

\begin{figure}[h!]
\hskip0.5cm
\begin{minipage}{5.5cm}
\includegraphics[width=5.5cm]{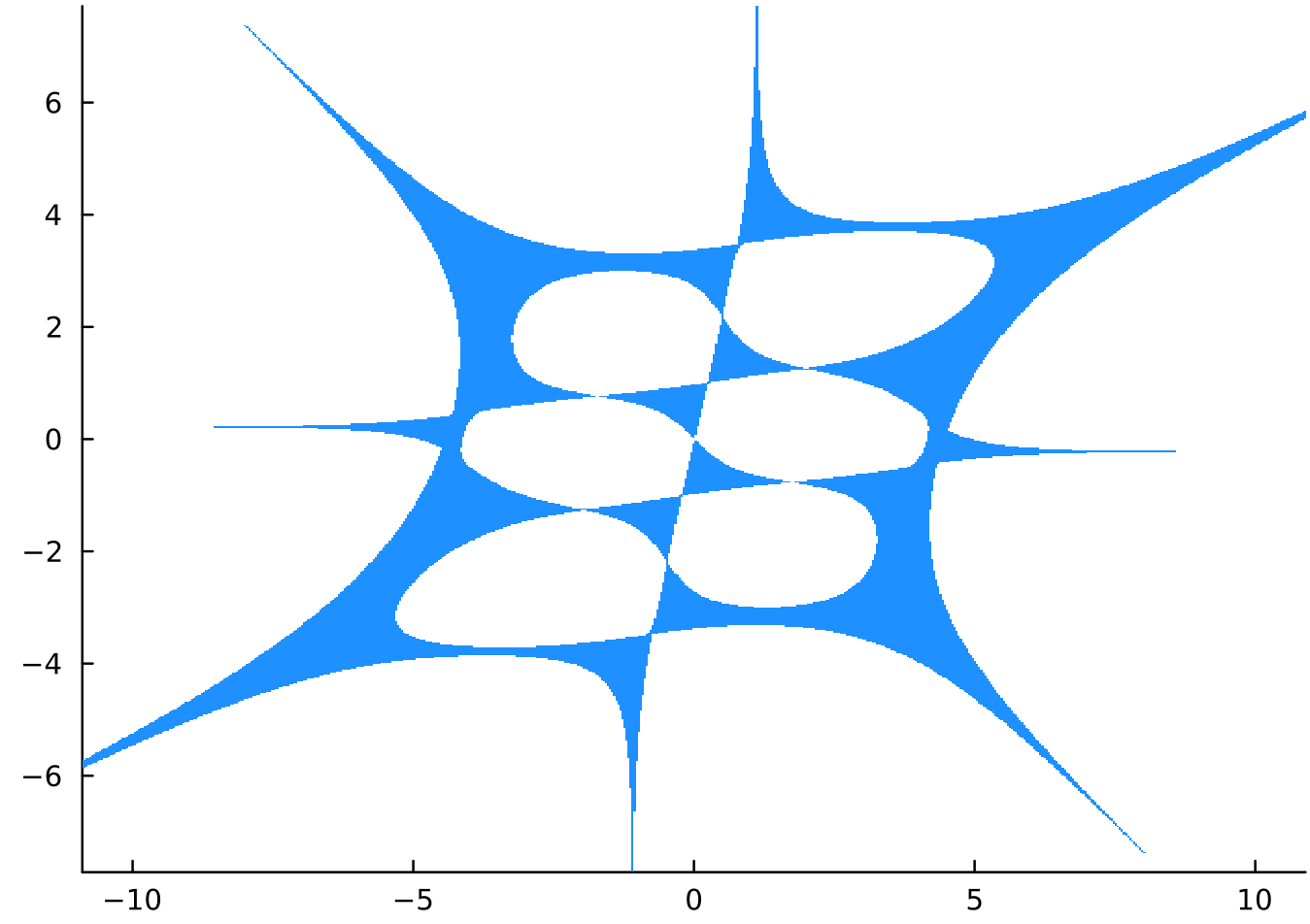}
\centering $(a)$
\end{minipage}
\hskip1cm
\begin{minipage}{5.5cm}
\includegraphics[width=5.5cm]{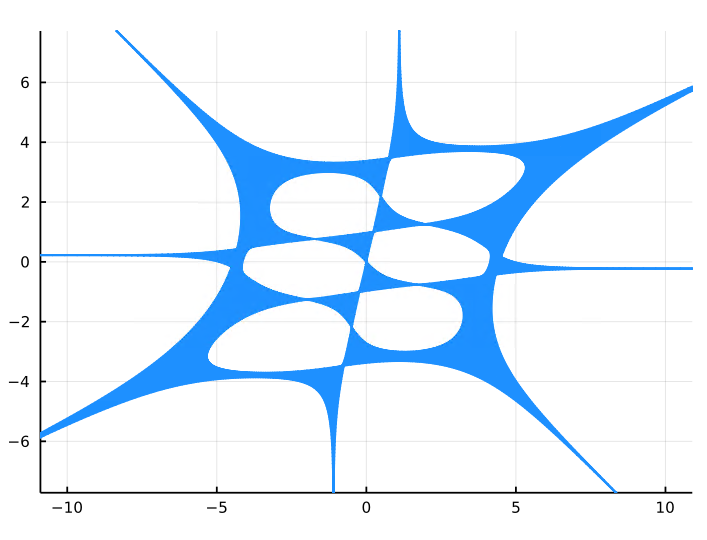}
\centering $(b)$
\end{minipage}
\caption{The amoeba of the polynomial~$p_3$ computed with the software [b] -- the ``greedy'' algorithm (a), the ``polygonal'' algorithm (b)}
\label{fig:p2_software_b}
\end{figure}

The result for the software [c] is shown in Figure~\ref{fig:p2_software_c} for two different values of the iterations parameter: 5 iterations, which is the minimal value such that the number of connected components of the amoeba complement in the picture coincides with the actual number, and 9 iterations, when the picture is smooth and the algorithm terminates much slower.

\begin{figure}[h!]
\hskip0.5cm
\begin{minipage}{5cm}
\includegraphics[width=5cm]{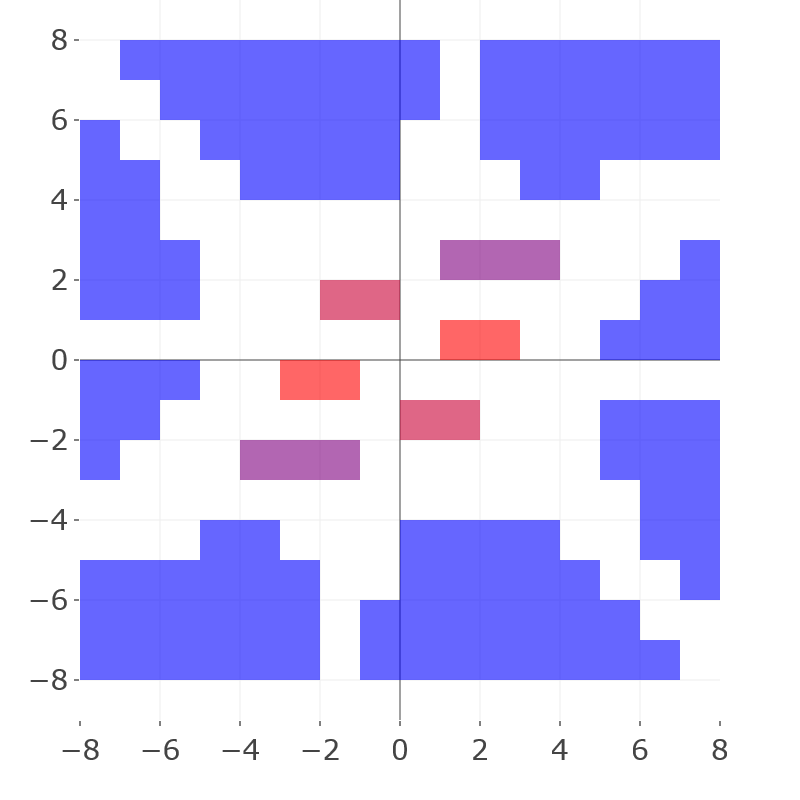}
\centering $(a)$
\end{minipage}
\hskip1cm
\begin{minipage}{5cm}
\includegraphics[width=5cm]{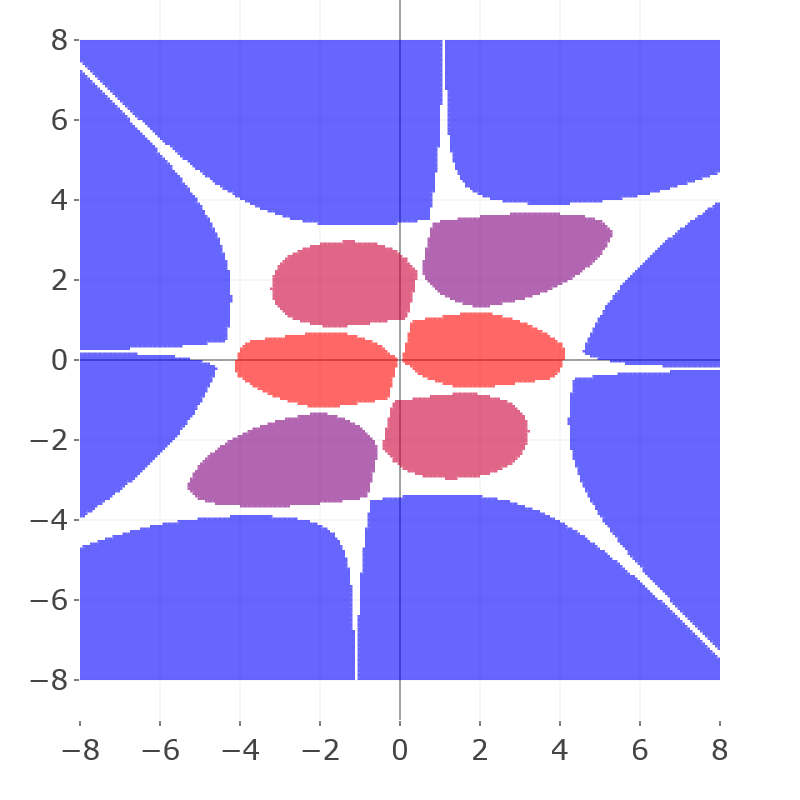}
\centering $(b)$
\end{minipage}

\caption{Connected components of the complement of the amoeba of the polynomial~$p_3$ computed with the software [c] -- 5~iterations (a), 9~iterations (b)}
\label{fig:p2_software_c}
\end{figure}
\end{example}

\subsubsection{Tests cases}

This subsection contains the description of the classes of polynomials that can be used to test the amoeba visualization software tools.

There are some classes of polynomials whose amoebas are of particular interest. Amoebas of optimal polynomials have the most complex structure so they are likely to be the most difficult to compute. Another aspect of complexity is a possibility for an amoeba complement to contain small bounded components (see Section~\ref{sec:small_components}).

It is also necessary to note one feature of any amoeba visualization software -- as the complexity of the depicted amoeba increases, so does the probability that the resulting image would contain various artifacts. These artifacts can look very different and, presumably, causes for their appearance also vary. In most cases, such artifacts can be easily detected, nevertheless, their presence in the image can make it difficult to analyze the structure of the amoeba.

\begin{example} \rm

Examples of image artifacts are presented in Figure~\ref{fig:artifacts}. Figure~\ref{fig:artifacts} (a) demonstrates the amoeba with missing lines in the center and underdrawn tentacles parallel to the y-axis, probably, due to its complexity. Another case, when there are extra lines in the picture, is shown in Figure~\ref{fig:artifacts}~(b).
\begin{figure}[h!]
\hskip0.5cm
\begin{minipage}{6cm}
\includegraphics[width=6cm]{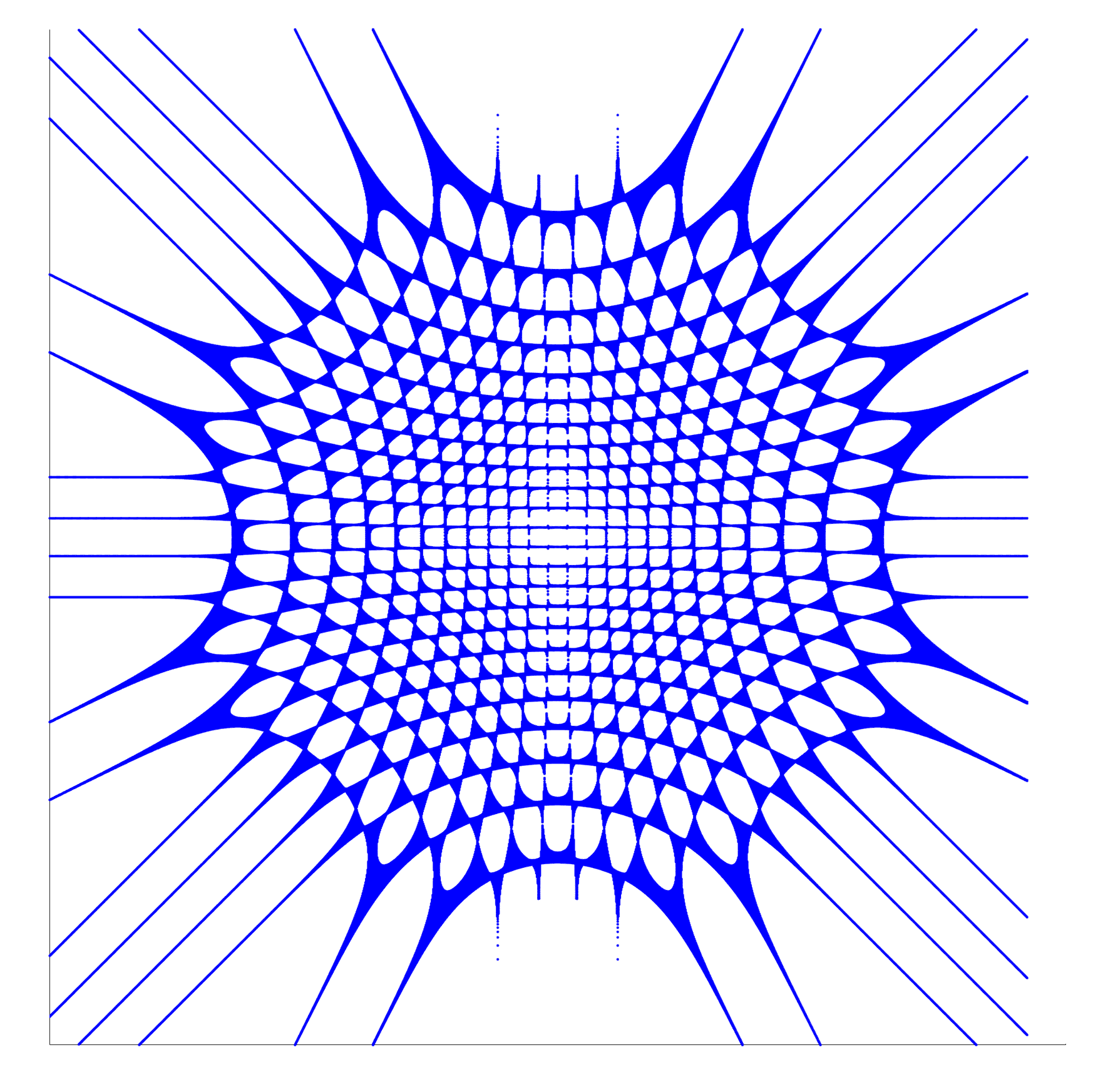}
\centering $(a)$
\end{minipage}
\hskip0.5cm
\begin{minipage}{6.25cm}
\includegraphics[width=6.25cm]{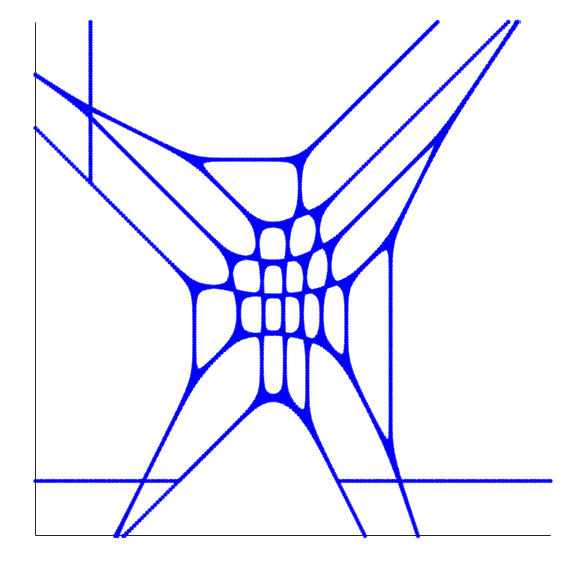}
\centering $(b)$
\end{minipage}

\caption{Image artifacts in the amoeba visualization -- missing lines in the picture (a), extra lines in the picture (b)}
\label{fig:artifacts}
\end{figure}
\label{ex:artifacts}
\end{example}

To test algorithms against polynomials whose amoebas have visual artifacts, it is proposed to use sets of random polynomials of a given degree, since the randomness of input polynomials leads to highly unpleasant amoebas with numerous tentacles.

Considering the above, the following sets of polynomials can be used for testing the software:

\begin{enumerate}
\item Optimal polynomials generated by the means of the algorithms presented in~\cite{Zhukov-Sadykov:2023}.
\item Polynomials whose amoeba complements contain small bounded components -- one can use the software [c] to generate such polynomials by using its amoeba evolution functionality. 
\item Sets of randomly generated polynomials.
\end{enumerate}

\subsection{Experimental conclusions}

This subsection contains the results of software testing on the test cases described above.

The software [b] seems to be inappropriate for the mass tests, since a lot of polynomials given as an input just lead to generating an exception. Possible reason for this is an update of the Julia packages combined with the lack of a support from the developer. The software [b] in the present state fits only for computing in some particular cases, since some input polynomials only generate an error and there were some cases, when the algorithm went into the infinite loop for no obvious reasons. For suitable polynomials the computation time is better than for the software [c].

The software [a] is faster than the software [c] in the case of optimal polynomials, but it has limitations on the degree of a polynomial, which becomes stricter with a growth of the number of monomials. For random polynomials there are a lot of issues (including image artifacts and even computation errors) with the result. 

The software [c] for the parameter values which ensure the best performance is slower than the software [a] and the computation time depends on the number of monomials, but it has lesser limitations on the parameters of input polynomials than other two packages.

\section{Conclusion}

The increase in computer performance significantly improves the capabilities of calculating the polynomial amoebas. The main complexity of this problem, however, lies more in the area of computational complexity of the algorithms.

The membership problem, though in general not solvable in polynomial time, has several different algorithmic solutions which are implemented in the software tools for the computation of polynomial amoebas. These solutions allow one to compute amoebas with a simple structure, but when computing amoebas of arbitrary polynomials, additional problems arise, such as the problem of identifying small components and appearance of artifacts (see Figure~\ref{fig:artifacts}).

The problem of small components is one of the keys to the Passare conjecture and it seems to have a high complexity by itself, since no software tools provide algorithms for its solution, though this problem is mentioned by the authors.~\cite{Timme-Master,Zhukov-Sadykov:2023}

Summing up the above, the answer to the main question formulated in the introduction is as follows. Though it is not clear at the moment whether the Passare conjecture is true or false, numerous computer experiments show that if a counterexample for it exists, then it is likely to be in higher dimensions. In order to find it, it is necessary to develop existing algorithms for constructing amoebas of polynomials having a high number of variables and their sections.

 \bibliographystyle{splncs04}
 \bibliography{amoebas}

\end{document}